\documentclass[12pt,twoside]{article}
\usepackage{graphicx}
\usepackage{amsfonts,amsbsy,amsthm,amssymb,amsmath,tikz}
\allowdisplaybreaks
\usepackage{cite}
\usepackage{graphics}
\def\bpsp{\begin{pspicture}}
\def\epsp{\end{pspicture}}

\newtheorem{theorem}{Theorem}[section]
\newtheorem{remark}[theorem]{Remark}
\newtheorem{example}[theorem]{Example}
\newtheorem{lemma}[theorem]{Lemma}
\newtheorem{corollary}[theorem]{Corollary}
\newtheorem{definition}[theorem]{Definition}
\newtheorem{proposition}[theorem]{Proposition}

\newtheorem{note}{Note}

\newtheorem{case}{Case}
\newtheorem{conjecture}{Conjecture}
\newtheorem{question}{Question}

\newcommand{\bea}{\begin{eqnarray}}
\newcommand{\eea}{\end{eqnarray}}
\newcommand{\beq}{\begin{eqnarray*}}
\newcommand{\eeq}{\end{eqnarray*}}

\def\m4{\mbox{\rm ~(mod $4$)}}

\def \bd{\begin{definition}}
\def \ed{\end{definition}}
\def \bqu{\begin{question}}
\def \equ{\end{question}}
\def \bcc{\begin{conjecture}}
\def \ecc{\end{conjecture}}
\def \bt{\begin{theorem}}
\def \et{\end{theorem}}
\def \bl{\begin{lemma}}
\def \el{\end{lemma}}
\def \bc{\begin{corollary}}
\def \ec{\end{corollary}}
\def \be{\begin{equation}}
\def \ee{\end{equation}}
\def \ben{\begin{enumerate}}
\def \een{\end{enumerate}}
\def \ba{\begin{array}}
\def \ea{\end{array}}
\def \bp{\begin{proposition}}
\def \ep{\end{proposition}}
\def \bx{\begin{example}}
\def \ex{\end{example}}
\def \br{\begin{remark}}
\def \er{\end{remark}}
\def \bdsc{\begin{description}}
\def \edsc{\end{description}}

\def \bn{\begin{case}}
\def \en{\end{case}}
\def \bnt{\begin{note}}
\def \ent{\end{note}}
\def\1{1\!\!1}

\def\mm2{\mbox{\rm ~(mod $2$)}}
\def\m4{\mbox{\rm ~(mod $4$)}}

\def\qed{\nolinebreak\hfill\rule{.2cm}{.2cm}\par\addvspace{.5cm}}

\def\m{\mu}

\def\1{\textbf{1}}
\def\0{\textbf{0}}

\begin{document}
\title{On distance Laplacian spread and Wiener index of a graph}
\author{ Saleem Khan$ ^{a} $, S. Pirzada$ ^{b} $\\
$^{a,b}${\em Department of Mathematics, University of Kashmir, Srinagar, Kashmir, India}\\
$ ^{a} $\texttt{khansaleem1727@gmail.com}; $^{b}$\texttt{pirzadasd@kashmiruniversity.ac.in}
}
\date{}

\pagestyle{myheadings} \markboth{Khan, Pirzada}{On distance Laplacian spread of graphs}
\maketitle
\vskip 5mm
\noindent{\footnotesize \bf Abstract.} Let $G$ be a simple connected simple graph of order $n$. The distance Laplacian matrix $D^{L}(G)$ is defined as $D^L(G)=Diag(Tr)-D(G)$, where $Diag(Tr)$ is the diagonal matrix of vertex transmissions and $D(G)$ is the distance matrix of $G$. The eigenvalues of $D^{L}(G)$ are the distance Laplacian eigenvalues of $G$ and are denoted by $\partial_{1}^{L}(G), \partial_{2}^{L}(G),\dots,\partial_{n}^{L}(G)$. The \textit{ distance Laplacian spread} $DLS(G)$ of a connected graph $G$ is the difference between largest and second smallest distance Laplacian eigenvalues, that is, $\partial_{1}^{L}(G)-\partial_{n-1}^{L}(G)$. We obtain bounds for $DLS(G)$ in terms of the Wiener index $W(G)$, order $n$ and the maximum transmission degree $Tr_{max}(G)$ of $G$ and characterize the extremal graphs. We obtain two lower bounds for $DLS(G)$, the first one in terms of the order, diameter and the Wiener index of the graph, and the second one in terms of the order, maximum degree and the independence number of the graph. For a connected $ k-partite$ graph $G$, $k\leq n-1$, with $n$ vertices having disconnected complement, we show that
 $ DLS(G)\geq \Big \lfloor \frac{n}{k}\Big \rfloor$
 with equality if and only if $G$ is a $complete ~ k-partite$ graph having cardinality of each independent class  same and $n \equiv 0 \pmod k$.

\vskip 3mm



\section{Introduction}\label{sec1}

Throughout this paper, we consider simple connected graphs. A simple connected graph  $G=(V,E)$ consists of the vertex set $V(G)=\{v_{1},v_{2},\ldots,v_{n}\}$ and the edge set  $E(G)$. The \textit{order} and \textit{size} of $G$ are $|V(G)|=n$ and  $|E(G)|=m$, respectively. The \textit{degree} of a vertex $v,$ denoted by $d_{G}(v)$ (we simply write by $d_v$) is the number of edges incident on $v$.  For other standard definitions, we refer to \cite{R2,R3}.\\
\indent Let $A$ be the adjacency matrix of $G$ and $D(G)=diag(d_1 ,d_2 ,\dots,d_n)$ the diagonal matrix of vertex degrees. The $Laplacian$ $matrix$ of $G$ is defined as $D(G)-A$. By the spectrum of $G$ we mean the spectrum of its adjacency matrix, and consists of the eigenvalues $\lambda_1 \geq \lambda_2 \geq \dots \geq \lambda_n$. The Laplacian spectrum of $G$ is the spectrum of its Laplacian matrix, and is denoted by $\mu_1 (G) \geq \mu_2 (G) \geq \dots \geq \mu_n (G) =0$. Laplacian eigenvalues are referred as $L-eigenvalues$ of $G$.\\
\indent In $G$, the \textit{distance} between two vertices $u,v\in V(G),$ denoted by $d_{uv}$, is defined as the length of a shortest path between $u$ and $v$. The \textit{diameter} of $G$, denoted by $diam(G)$, is the maximum distance between any two vertices of $G.$ The \textit{distance matrix} of $G$ is denoted by $D(G)$ and is defined as $D(G)=(d_{uv})_{u,v\in V(G)}$.
The \textit{transmission} $Tr_{G}(v)$
(or briefly $Tr(v)$ if graph $G$ is understood) of a vertex $v$ is defined as the sum of the distances from $v$ to all other vertices in $G$, that is, $Tr_{G}(v)=\sum\limits_{u\in V(G)}d_{uv}.$ A graph $ G $ is said to be $k$-\textit{transmission regular} if $ Tr_{G}(v)=k,$ for each $ v\in V(G)$.   The \textit{Wiener index}  (also called \textit{transmission})  of a graph $G$,  denoted by $W(G)$,  is the sum of distances between all unordered pairs of vertices in $G$.  Clearly, $W(G)=\frac{1}{2}\displaystyle\sum_{u,v\in V(G)}d_{uv}=\frac{1}{2}\displaystyle\sum_{v\in V(G)}Tr_{G}(v). $ For any vertex $v_i\in V(G)$, the transmission $Tr_G(v_i)$ is also called the \textit{transmission degree}.

Let $Tr(G)=diag (Tr_1,Tr_2,\ldots,Tr_n) $ be the diagonal matrix of vertex transmissions of $G$. Aouchiche and Hansen \cite{R1} introduced the Laplacian for the distance matrix of a connected graph. The matrix $D^L(G)=Tr(G)-D(G)$ (or simply $D^{L}$) is called the \textit{distance Laplacian matrix} of $G$. The eigenvalues of $D^{L}(G)$ are called the distance Laplacian eigenvalues of $G$ and are referred as $D^L-eigenvalues$ of $G$. Since $ D^L(G) $ is a real symmetric positive semi-definite matrix, we take the distance Laplacian eigenvalues as $0=\partial_{n}^{L}(G)\leq \partial_{n-1}^{L}(G)\leq \dots\leq \partial_{1}^{L}(G)$. If $\partial_{i}^{L}(G)$ is repeated $p$ times, then we say multiplicity of $\partial_{i}^{L}(G)$ is $p$ and we write $m(\partial_{i}^{L}(G))=p$. \\
\indent As usual, we denote by $K_n$, $S_n$ and $K_{n-2}\vee \overline{K_2}=K_{n}-e$  the complete graph, the star and the complement of an edge, respectively, all on $n$ vertices.\\
\indent The spread of a graph $G$ is defined as $(G)=\lambda_1 -\lambda_n$, where $\lambda_1$ and $\lambda_n$ are respectively the largest and the smallest eigenvalues of the adjacency matrix of $G$. The Laplacian spread of a graph $G$ is defined as
$LS(G)=\mu_1 -\mu_{n-1}$, where $\mu_1$ and $\mu_{n-1}$ are respectively the largest and second smallest Laplacian eigenvalues of $G$.  More on $ LS(G)$ can be found in \cite{R10,R11,R12}.

The distance Laplacian spread of a graph $G$, denoted   by $ DLS(G)$, is defined as
$$DLS(G)=\partial_{1}^{L}(G)-\partial_{n-1}^{L}(G),$$
where $\partial_{1}^{L}(G)$ and $\partial_{n-1}^{L}(G)$ are respectively the largest and second smallest distance Laplacian eigenvalues of $G$.

The rest of the paper is organized as follows. In Section 2, we state some results which will be used in the sequel. In Section 3, we obtain lower and upper bounds for $DLS(G)$ in terms of the Wiener index $W(G)$, order $n$ and the maximum transmission degree $Tr_{max}(G)$ of $G$ and characterize the extremal graphs. In Section 4, we obtain two lower bounds for $DLS(G)$, the first one in terms of the order, diameter and the Wiener index of the graph, and the second one in terms of the order, maximum degree and the independence number of the graph. For a connected $ k-partite$ graph $G$, $k\leq n-1$, with $n$ vertices having disconnected complement, we show that
 $ DLS(G)\geq \Big \lfloor \frac{n}{k}\Big \rfloor$
 with equality if and only if $G$ is a $complete ~ k-partite$ graph having cardinality of each independent class  same and $n \equiv 0 \pmod k$. Also, for a connected graph $G$ with given chromatic number $\chi$, $\frac{n}{2}\leq \chi \leq n-1$, and complement of $G$ being disconnected, we prove that $DLS(G)\geq 2$ with equality if and only if $G\cong K_{\underbrace{{2,2,\dots,2}}_{n-\chi},\underbrace{{1,1,\dots,1}}_{2\chi -n}}$.

\section{Preliminaries }

In this section, we state some lemmas which will be used to prove our main results.

\begin{lemma}\label{L1}\emph{\cite{R1}} Let $G$ be a connected graph on $n$ vertices. Then $\partial^L_{n-1}(G)=n$ if and only if $\overline {G}$ is disconnected. Furthermore, the multiplicity of $n$ as a distance Laplacian eigenvalue is one less than the number of connected components of $\overline {G}$.
\end{lemma}
\begin{lemma}\label{L2}\emph{\cite{R4}} Let $G$ be a connected graph with $n$ vertices having $\partial^L_{1}(G)$ as the distance Laplacian spectral radius. Then $\partial^L_{1}(G)\geq max\{Tr(v):v\in V(G)\}+1 $ with equality if and only if $G\cong K_n$.
\end{lemma}
\begin{lemma}\label{L3}\emph {\cite{R5}} Let $G$ be a connected graph such that $D^L$ has an eigenvalue with multiplicity $n-2$. Then, exactly one of the following condition holds.\\
(i) $m(\partial^L_{1})=n-2$ and $G\cong S_n$ or $G\cong K_{p,p}$.\\
(ii) $m(\partial^L_{n-1})=n-2$ and $G\cong K_{n-2}\vee \overline K_{2}$.
\end{lemma}
\begin{lemma}\label{L4} \emph{\cite{R8}}  Let $t_{1},t_{2},\dots,t_{k}$ and n be integers such that $t_{1}+t_{2}+\dots+t_{k}=n$ and $t_{i}\geq 1$ for $i=1,2,\dots,k$. Let $p=|\{i:t_{i}\geq 2\}|$. The distance Laplacian spectrum of the complete $k-partite$ graph $K_{t_{1},t_{2},\dots,t_{k}}$ is$ \Big({(n+t_{1})}^{(t_{1}-1)},\dots,{(n+t_{p})}^{(t_{p}-1)},n^{(k-1)},0\Big)$.
\end{lemma}
\begin{lemma}\label{L5} \emph{\cite{R1}}  Let $G$ be a connected graph with $n$ vertices and $m$ edges, where $m\geq n$. Let $G^*$ be the connected graph obtained from $G$ by deleting an edge. Let $\partial^L_1 \geq \partial^L_2 \geq ...\geq \partial^L_n$ and ${\partial^*_1}^L \geq {\partial^*_2}^L \geq ...\geq {\partial^*_n}^L$ be the spectrum of $G$ and $G^*$, respectively. Then ${\partial^*_i}^L \geq \partial^L_i $ for all $i=1,\dots,n$.
\end{lemma}

\section{Bounds for the distance Laplacian spread of a graph}

For a graph $G$ with $n$ vertices, let $Tr_{max}(G)=max\{Tr(v):v\in V(G)\}$ and $Tr_{min}(G)=min\{Tr(v):v\in V(G)\}$. Whenever the graph $G$ is understood, we will write $Tr_{max}$ and $Tr_{min}$ in place of $Tr_{max}(G)$ and $Tr_{min}(G)$, respectively.  For every connected graph $G$, we know that 0 is always a distance Laplacian eigenvalue with multiplicity one. Therefore, from the definitions, $2W=\partial^L_{1}(G)+\partial^L_{2}(G)+\dots+\partial^L_{n-1}(G)$. Also, $Tr_{max}\geq \frac{2W}{n}$ and  $Tr_{min}\leq \frac{2W}{n}$, where $\frac{2W}{n}$ is the average transmission degree.

Clearly ${(\partial^L_{1}(G))}^2, {(\partial^L_{2}(G))}^2 ,\dots,{(\partial^L_{n-1}(G))}^2$ are the non-zero distance Laplacian eigenvalues of ${(D^L (G))}^2$. Using the fact that the trace of a matrix is equal to sum of its  eigenvalues, we observe that
\begin{equation}
 \sum_{i=1}^{n-1}{{(\partial^L_{i}(G))}^2=\sum\limits_{i}Tr^2(i)+2\sum\limits_{1\leq i<j\leq n}d^2_{ij}}
\end{equation}
\indent We recall the following well-known result from matrix theory.
\begin{lemma}\label{IL1}\emph {\cite{R6}} Let $M=(m_{ij})$ be a complex $n\times n$-matrix having $l_1 ,l_2 ,\dots,l_p$ as its distinct eigenvalues. Then
$$\{l_1, l_2,\dots,l_p\}\subset \bigcup\limits_{i=1}^{n}\Big \{z:|z-m_{ii}|\leq \sum\limits_{j\neq i}|m_{ij}|\Big\}.$$
\end{lemma}
\indent By application of Lemma \ref{IL1} for the distance Laplacian matrix of a graph $G$ with $n$ vertices, we easily get the inequality
\begin{equation}
\partial^L_{1}(G)\leq 2Tr_{max}
\end{equation}

\indent Now, we obtain the lower and upper bounds for the distance Laplacian spread $DLS(G)$ of a graph $G$ in terms of the Wiener index $W(G)$, order $n$ and the maximum transmission degree $Tr_{max}(G)$ of the graph $G$.

\begin{theorem}\label{T1}
Let $G$ be a connected graph with $n$ vertices having Wiener index $W(G)$. Then
$$\frac{(n-1)(Tr_{max}+1)-2W(G)}{n-2}\leq DLS(G)\leq 2(n-1)Tr_{max}-2W(G)$$
Equality holds in the left if and only if $G\cong K_n$.
\end{theorem}
\noindent{\bf Proof.} Let $\partial_{1}^{L}(G),\partial_{2}^{L}(G),\dots,\partial_{n}^{L}(G)$ be the $D^L-eigenvalues$ of $G$. As $\partial_{n}^{L}(G)=0$ is always a $D^L-eigenvalue$ of $G$ with multiplicity one, therefore, we have
$2W(G)=\partial_{1}^{L}(G)+\partial_{2}^{L}(G)+\dots+\partial_{n-1}^{L}(G)$, which gives $2W(G)\geq \partial_{1}^{L}(G)+(n-2)\partial_{n-1}^{L}(G)$, further implies that $\partial_{n-1}^{L}(G)\leq \frac{2W(G)-\partial_{1}^{L}(G)}{n-2}$.
Hence,\\
\begin{align*}
DLS(G)&=\partial_{1}^{L}(G)-\partial_{n-1}^{L}(G)
\geq \partial_{1}^{L}(G)-\frac{2W(G)-\partial_{1}^{L}(G)}{n-2}\\&
 =\frac{(n-1)\partial^L_{1}(G)-2W(G)}{n-2}
\geq \frac{(n-1)(Tr_{max}+1)-2W(G)}{n-2} ~~~~~~~(using ~  Lemma ~ \ref{L2})
\end{align*}
This proves the left hand side inequality.

Clearly equality holds in above if and only if $\partial_{1}^{L}(G)\geq \partial_{2}^{L}(G)=\dots =\partial_{n-1}^{L}(G)$ and $G\cong K_n$. We consider the following two cases.\\
\noindent{\bf Case 1.} Let  $\partial_{1}^{L}(G)= \partial_{2}^{L}(G)=\dots =\partial_{n-1}^{L}(G)$. Thus $m(\partial_{1}^{L}(G))=n-1$ which shows that $G\cong K_n$.\\
\noindent{\bf Case 2.} Let  $\partial_{1}^{L}(G)> \partial_{2}^{L}(G)=\dots =\partial_{n-1}^{L}(G)$. Using Lemma \ref{L3}, we get $G\cong K_{n-2}\vee \overline K_{2}$.\\
\indent Thus equality holds in the left hand side of the inequality if and only if $G\cong K_n$.

Now to prove the right hand inequality, we consider $2W(G)=\partial_{1}^{L}(G)+\partial_{2}^{L}(G)+\dots+\partial_{n-1}^{L}(G)$, which gives $2W(G)\leq (n-2)\partial_{1}^{L}(G)+\partial_{n-1}^{L}(G)$, further implies that $\partial_{n-1}^{L}(G)\geq 2W(G)-(n-2)\partial_{1}^{L}(G).$ Therefore,
\begin{align*}
DLS(G)&=\partial_{1}^{L}(G)-\partial_{n-1}^{L}(G) \leq \partial_{1}^{L}(G)-2W(G)+(n-2)\partial_{1}^{L}(G)\\
&  =(n-1)\partial^L_{1}(G)-2W(G) \leq 2(n-1)Tr_{max}-2W(G) ~~~~~~~~ using ~ inequality ~ (3.2)
\end{align*}
proving the right hand inequality. \qed

\begin{corollary}
Let $G$ be a connected graph with $n$ vertices having Wiener index $W(G)$. Then
\begin{equation}
DLS(G)\geq \frac{n^2 -n-2W(G)}{n^2 -2n}
\end{equation}
with equality if and only if $G\cong K_n$.
\end{corollary}
\noindent{\bf Proof.} Using the fact that the maximum transmission degree is always greater or equal to average transmission degree, that is, $Tr_{max}\geq \frac{2W(G)}{n}$ in the left hand inequality of Theorem \ref{T1}, we get
\begin{equation*}
DLS(G)\geq \frac{(n-1)(\frac{2W(G)}{n}+1)-2W(G)}{n-2},
\end{equation*}
or
\begin{equation*}
DLS(G)\geq \frac{n^2 -n-2W(G)}{n^2 -2n}.
\end{equation*}
Thus inequality (3.3) follows.\\
We note that the average transmission degree is  equal to maximum transmission degree for a complete graph $K_n$, therefore we observe that equality holds in (3.3) if and only if $G\cong K_n$. \qed

\begin{theorem}\label{T2}
Let $G$ be a connected graph with $n$ vertices having Wiener index $W(G)$.  If the complement of $G$ is disconnected, then\\
$$1+Tr_{max}-n\leq DLS(G)\leq 2Tr_{max}-n$$
 Equality holds in the left if and only if $G\cong K_n$.
\end{theorem}
\noindent{\bf Proof.} Assume that $\overline{G}$ is disconnected. Using Lemma \ref{L1}, we get $\partial_{n-1}^{L}(G)=n$ is the $D^L -eigenvalue$ of $G$ with some multiplicity greater or equal to one depending on the number of components in $\overline{G}$. Thus, using Lemma \ref{L2} and inequality (3.2), we get
\begin{equation*}
1+Tr_{max}-n\leq DLS(G)=\partial_{1}^{L}(G)-\partial_{n-1}^{L}(G)\leq 2Tr_{max}-n
\end{equation*}
with equality holding in the left if and only if $G\cong K_n$. \qed

\begin{corollary}
Let $G$ be a connected graph with $n$ vertices having Wiener index $W(G)$.  If the complement of $G$ is disconnected, then
$ DLS(G)\geq \frac{2W(G)-n(n-1)}{n}$ with equality if and only if $G\cong K_n$.
\end{corollary}
\noindent{\bf Proof.} Using the fact that $Tr_{max}\geq \frac{2W(G)}{n}$ and the observation that the average transmission degree is equal to the maximum transmission degree for a complete graph in \ref{T2}, we have
\begin{equation*} DLS(G)\geq 1+Tr_{max}-n \geq 1+ \frac{2W(G)}{n}-n= \frac{2W(G)-n(n-1)}{n}
\end{equation*}
with equality if and only if $G\cong K_n$. \qed

The following result gives a lower bound for the distance Laplacian spread of the graph $G$ in terms of $Tr_{max}$ and order  $n$.
\begin{theorem}\label{T3}
Let $G$ be a connected graph with $n$ vertices. Then
\begin{equation*}
DLS(G)\geq 1+Tr_{max}-\sqrt{\frac{R_1 -{(1+Tr_{max})}^2}{n-2}},
\end{equation*}
where $R_{1}=\sum\limits_{i}Tr^2(i)+2\sum\limits_{1\leq i<j\leq n}d^2_{ij}$. Equality holds  if and only if $G\cong K_n$.
\end{theorem}
\noindent{\bf Proof.} From equation (3.1), we have
$$\sum_{i=1}^{n-1}{{(\partial^L_{i}(G))}^2=\sum\limits_{i}Tr^2(i)+2\sum\limits_{1\leq i<j\leq n}d^2_{ij}}=R_1$$
Clearly, $R_1=\sum_{i=1}^{n-1}{(\partial^L_{i}(G))}^2\geq {(\partial^L_{1}(G))}^2 +(n-2){(\partial^L_{n-1}(G))}^2$,
 or $\partial^L_{n-1}(G)\leq \sqrt{\frac{R_1 -{(\partial^L_{1}(G))}^2}{n-2}}$, that is, $ \partial^L_{n-1}(G)\leq \sqrt{\frac{R_1 -{(\partial^L_{1}(G))}^2}{n-2}}.$ Using this inequality for $\partial^L_{n-1}(G)$, we have
 \begin{equation}
 DLS(G)=\partial^L_{1}(G)-\partial^L_{n-1}(G)\geq \partial^L_{1}(G)-\sqrt{\frac{R_1 -{(\partial^L_{1}(G))}^2}{n-2}}
 \end{equation}
We note that $f(y)=y-\sqrt{\frac{R_1 -{y}^2}{n-2}}$ is an increasing function. Using this fact and Lemma \ref{L2} in inequality (3.4), we get
 \begin{equation}
  DLS(G)\geq 1+Tr_{max}-\sqrt{\frac{R_1 -{(1+Tr_{max})}^2}{n-2}}.
  \end{equation}
  This establishes the inequality.\\

Clearly equality holds in (3.5) if and only if $\partial_{1}^{L}(G)\geq \partial_{2}^{L}(G)=\dots =\partial_{n-1}^{L}(G)$ and $G\cong K_n$. Using the fact that complete graph is the only graph with two distinct distance Laplacian eigenvalues and Lemma \ref{L2}, it follows that equality holds in (3.4) if and only if  $G\cong K_n$. \qed

\begin{corollary}
Let $G$ be a connected graph with $n$ vertices. Then
\begin{equation*} DLS(G)\geq \frac{n+2W(G)}{n}-\frac{\sqrt{n^2 R_1 +{(n+2W(G))}^2}}{n\sqrt{n-2}}.
 \end{equation*}
Equality holds  if and only if $G\cong K_n$.
\end{corollary}
\noindent{\bf Proof.} Using the arguments as in Corollary 3.5, the result is established. \qed

\section{Distance Laplacian spread of graphs with given diameter and independence number}

We start this section with the following theorem.

\begin{theorem}\label{T5}
(Cauchy Interlacing Theorem). Let $M$ be a real symmetric matrix of order $n$, and let $A$ be a principal submatrix of $M$ with order $s\leq n$. Then $$\lambda_i (M)\geq \lambda_i (A) \geq \lambda_{i+n-s} (A)\hspace{1cm}(1\leq i\leq s).$$
\end{theorem}

Now, we obtain a lower bound for $DLS(G)$ in terms of order $n$, diameter $d$ and Wiener index $W(G)$.

\begin{theorem}\label{T7}
Let $G$ be a connected graph with $n$ vertices having diameter $d$ and Wiener index $W(G)$. Then
\begin{equation*}
DLS(G)\geq \frac{2n+d^2 -2d+1}{2}-\frac{2W(G)}{n-1}.
\end{equation*}
\end{theorem}
\noindent{\bf Proof.} As the diameter of $G$ is $d$, we assume that $v_1 ,v_2 ,\dots,v_{d+1}$ be the vertices of the diametral path in that order. It is easy to see that
\begin{eqnarray}
\begin{split}
Tr(v_1)&\geq 1+2+\dots+d+n-d-1\\
&=\frac{d(d+1)}{2}+n-d-1\\
&=\frac{2n+d^2 -d-2}{2}
\end{split}
\end{eqnarray}
and
\begin{eqnarray}
\begin{split}
Tr(v_2)&\geq 1+ 1+2+\dots+d-1+n-d-1\\
&=\frac{(d-1)d}{2}+n-d\\
&=\frac{2n+d^2 -3d}{2}
\end{split}
\end{eqnarray}
Also, the principal submatrix of $D^L (G)$, say $L$, corresponding to vertices $v_1$ and $v_2$ is given by
\begin{equation*}
L=
\begin{pmatrix}
Tr(v_1) & -1 \\
-1 & Tr(v_2) \\
\end{pmatrix}.
\end{equation*}
The characteristic polynomial of $L$ is given by
\begin{equation*}
f(y)=y^2 -(Tr(v_1)+Tr(v_2))y+Tr(v_1)Tr(v_2)-1.
\end{equation*}
Let $y_1$ and $y_2$ with $y_1\geq y_2$, be the roots of $f(y)$. Then, clearly
\begin{equation*}
y_1 =\frac{Tr(v_1)+Tr(v_2)+\sqrt{{(Tr(v_1)-Tr(v_2))}^2 +4}}{2}
\end{equation*}
so that
\begin{equation*}
y_1 \geq \frac{Tr(v_1)+Tr(v_2)+2}{2}.
\end{equation*}
Using Theorem \ref{T5}, we have
\begin{equation*}
\partial^L_1 (G)\geq y_1\geq \frac{Tr(v_1)+Tr(v_2)+2}{2}.
\end{equation*}
With the help of inequalities (4.6) and (4.7), we get
\begin{equation}
\partial^L_1 (G)\geq \frac{Tr(v_1)+Tr(v_2)+2}{2}\geq \frac{2n+d^2 -2d+1}{2}.
\end{equation}
Also, we have the following upper bound for $\partial^L_{n-1} (G)$ in terms of $W(G)$ and order $n$.
\begin{equation}
\partial^L_{n-1} (G)\leq \frac{\sum\limits_{i=1}^{n-1} \partial^L_i (G)} {n-1}=\frac{2W(G)}{n-1}
\end{equation}
Using inequalities (4.8) and (4.9), we get
$$DLS(G)=\partial^L_1 (G)-\partial^L_{n-1} (G)\geq \frac{2n+d^2 -2d+1}{2}-\frac{2W(G)}{n-1}$$ \qed

We require the following lemmas to prove the next result.

\begin{lemma}\label{IL2}\emph{\cite{R1}} Let $G$ be a connected graph on $n$ vertices with $diam(G)\leq 2$. Let $\mu_1 (G) \geq \mu_2 (G)\geq \dots \geq \mu_n (G)=0$ be the Laplacian spectrum of $G$. Then the distance Laplacian spectrum of $G$ is  $2n-\mu_{n-1} (G) \geq 2n- \mu_{n-2} (G)\geq \dots \geq 2n-\mu_1 (G)>\partial^L_n (G)=0$. Moreover, for every $i\in \{1,2,\dots,n-1\}$ the eigenspaces corresponding to $\mu_i (G)$ and $2n-\mu_i (G)$ are same.
\end{lemma}
\begin{lemma}\label{IL3}\emph{\cite{R9}}
Let $G$ be a graph with at least one edge and maximum vertex degree $\triangle (G)$. Then $\mu_1 (G) \geq 1+\triangle (G)$ with equality for a connected graph  if and only if $\triangle (G)=n-1$.
\end{lemma}

In a graph $G$, an independent set is a subset $S$ of the vertex set $V(G)$ if no two vertices of $S$ are adjacent. The $ independence ~ number$ of $G$, denoted by $\alpha (G)$, is defined as
$\alpha (G)=max\{|S|:S ~ is ~ an ~ independent ~ set ~ of ~ G\}$.

\begin{theorem}
Let $G$ be a connected graph with $n$ vertices having diameter $d=2$, maximum vertex degree $\triangle (G)$ and independence number $\alpha (G)$. Then
\begin{equation*}DLS(G)\geq \triangle (G)+ \alpha (G)+1 -n.
\end{equation*}
\end{theorem}
\noindent{\bf Proof.} Clearly, $G\ncong K_n$ and so $\alpha (G)\geq 2$. Let $\{v_1 ,v_2 ,\dots,v_{\alpha}\}$ be the maximal independent set. We have
\begin{eqnarray}
Tr(v_1 )\geq 2(\alpha -1) +n-\alpha=n+\alpha -2
\end{eqnarray}
Similarly, we have
\begin{equation}
Tr(v_2 )\geq n+\alpha -2
\end{equation}
The principal submatrix corresponding to the vertices $v_1$ and $v_2$ is given by
\begin{equation*}
P=
\begin{pmatrix}
Tr(v_1) & -2 \\
-2 & Tr(v_2) \\
\end{pmatrix}.
\end{equation*}
The characteristic polynomial corresponding to $P$ is given by
\begin{equation*}
f(x)=x^2 -(Tr(v_1 ) +Tr(v_2 ))x+Tr(v_1 )Tr(v_2 )-4.
\end{equation*}
Let $x_1$ and $x_2$ with $x_1 \geq x_2$, be the roots of $f(x)$. Then, we have
\begin{equation*}
x_1 =\frac{Tr(v_1)+Tr(v_2)+\sqrt{{(Tr(v_1)-Tr(v_2))}^2 +16}}{2},
\end{equation*}
so that
\begin{equation*}
x_1 \geq \frac{Tr(v_1)+Tr(v_2)+4}{2}.
\end{equation*}
Using Theorem \ref{T5}, we get
\begin{equation}
\partial^L_1 (G)\geq x_1\geq \frac{Tr(v_1)+Tr(v_2)+4}{2}.
\end{equation}
Now, by using  inequalities (4.10) and (4.11) in inequality (4.12), we get
\begin{equation}
\partial^L_1 (G)\geq n+\alpha.
\end{equation}
Using Lemmas \ref{IL2} and \ref{IL3}, we have
\begin{equation}
\partial^L_{n-1} (G)=2n-\mu_1(G) \leq 2n-1 -\triangle (G)
\end{equation}
Hence, from inequalities (4.13) and (4.14), we obtain
\begin{equation*}
DLS(G)=\partial^L_1 (G)-\partial^L_{n-1}\geq n+\alpha -2n+1+\triangle (G),
\end{equation*}
or
\begin{equation*}
DLS(G)\geq \triangle (G)+ \alpha (G)+1-n.
\end{equation*} \qed

For an integer $k\geq 2$, a $k-partite$ graph is a graph whose vertices can be partitioned into $k$ different independent classes such that every edge of the graph has ends in different independent classes. A $complete ~ k-partite$ graph is a $k-partite$ graph in which there is an edge between every pair of vertices from different independent classes. $Complete ~ k-partite$ graph is  usually denoted by $K_{l_1 ,l_2 ,\dots,l_k}$, where $k$ is the number of independent classes, $l_i$`s are the cardinalities of independent classes, $(i=1,2,\dots,k)$  and $l_1 +l_2 +...+l_k =n$. We will assume that $l_1 \geq l_2 \geq \dots \geq l_s$. If $k=2$, it is called a complete bipartite graph. As usual we denote a complete bipartite graph by $K_{n-r,~r}$, $1\leq r\leq n-1$. \\
\indent Also, for $x$ real, $\lfloor x\rfloor$ denotes the floor function which is equal to the greatest integer less than or equal to $x$.

The following theorem gives a lower bound for $DLS(G)$, when $G$ is a $ k-partite$ graph having disconnected complement.

\begin{theorem}\label{T6}
 Let $G$ be a connected $ k-partite$ graph, $k\leq n-1$, with $n$ vertices having disconnected complement. Then
 \begin{equation}
 DLS(G)\geq \Big \lfloor \frac{n}{k}\Big \rfloor
 \end{equation}
 with equality if and only if $G$ is $complete ~ k-partite$ graph having cardinality of each independent class  same with $n \equiv 0 \pmod k$.
\end{theorem}
\noindent{\bf Proof.} As $G$ is given to be a $k-partite$ graph, it is a spanning subgraph of $complete ~ k-partite$ graph $ K_{l_1 ,l_2 ,\dots,l_k}$, where $l_i$, $i=1,2,\dots,k$, is the cardinality of the $ith$ independent class and $l_1 +l_2 +...+l_k =n$ . Let us order the cardinalities of independent classes as $l_1  \geq l_2 \geq \dots \geq l_k$. Using Lemmas \ref{L5} and  \ref{L4}, we get
\begin{equation*}
\partial^L_1 (G)\geq \partial^L_1 ( K_{l_1 ,l_2 ,\dots,l_k})=n+l_1.
\end{equation*}
Since $l_1$ is the largest among the cardinalities of the independent classes, it is at least equal to average cardinality. Thus from above, we have
\begin{equation}
\partial^L_1 (G)\geq n+l_1 \geq n+\frac{n}{k}\geq n+\Big \lfloor \frac{n}{k}\Big \rfloor
\end{equation}
Also, $\overline{G}$ is disconnected, so that by Lemma \ref{L1}, $\partial^L_{n-1} (G)=n$ is an eigenvalue of $D^L (G)$. Therefore,
\begin{equation*}
DLS(G)=\partial^L_1 (G)-\partial^L_{n-1} (G)\geq n+\Big \lfloor \frac{n}{k}\Big \rfloor -n=\Big \lfloor \frac{n}{k}\Big \rfloor.
\end{equation*}
Thus the inequality is established.

Now, if $G\cong   K_{l_1 ,l_2 ,\dots,l_k}$ satisfying $l_1 =l_2 =\dots=l_k    $ and $n \equiv 0 \pmod k$, then using Lemma \ref{L4}, we have
\begin{equation*}
\partial^L_1 (G)= n+l_1 = n+\frac{n}{k}= n+\Big \lfloor \frac{n}{k}\Big \rfloor
\end{equation*}
so that
\begin{equation*}
DLS(G)=\partial^L_1 (G)-\partial^L_{n-1} (G)= n+\Big \lfloor \frac{n}{k}\Big \rfloor -n=\Big \lfloor \frac{n}{k}\Big \rfloor
\end{equation*}
and the equality holds.

To complete the proof, it suffices to show that  if the given $k-partite$  graph is different from the one mentioned in the statement of theorem, then the inequality (4.15) is strict. We have following cases to consider.

\noindent{\bf Case 1.} Assume that all $l_i$`s, $i=1,2,\dots,k$, are not equal, so that $l_1>\frac{n}{k} $. Thus, from inequality (4.16), $\partial^L_1 (G)>n+\Big \lfloor \frac{n}{k}\Big \rfloor$ which further shows that  $DLS(G)>\Big \lfloor \frac{n}{k}\Big \rfloor$.

\noindent{\bf Case 2.} Let $n \not\equiv 0 \pmod k$. Then, from inequality (4.16), we get $\partial^L_1 (G)\geq n+\frac{n}{k}>n+\Big \lfloor \frac{n}{k}\Big \rfloor$ which  shows that  $DLS(G)>\Big \lfloor \frac{n}{k}\Big \rfloor$.

\noindent{\bf Case 3.} Now, for the last case, we assume that $k\geq 3$ because of the fact that a complete bipartite graph is the only $2-partite$ graph having disconnected complement. Let $G$ be isomorphic to a proper spanning subgraph of $complete ~ k-partite$ graph $ K_{l_1 ,l_2 ,\dots,l_k}$, $k\geq 3$,  having disconnected complement with all $l_i$`s, $i=1,2,\dots,k$, being equal, that is, $l_1 =l_2 =\dots=l_k$ and  $n \equiv 0 \pmod k$. Using Lemma \ref{L5}, it suffices to prove that inequality (4.15) is strict for $G$ isomorphic to $ K_{l_1 ,l_2 ,\dots,l_k}-e  $, where $e$ is an edge between the vertices, say $u$ and $v$ of $ K_{l_1 ,l_2 ,\dots,l_k}  $. Clearly, $Tr(u)=Tr(v)= 2(l_1 -1)+2+n-l_1 -1=n+l_1 -1$.

The principal submatrix of $D^L(G)$ corresponding to vertices $u$ and $v$ is given by
\begin{equation*}
Q=
\begin{pmatrix}
Tr(u) & -2 \\
-2 & Tr(v) \\
\end{pmatrix}
=
\begin{pmatrix}
n+l_1 -1 & -2 \\
-2 & n+l_1 -1 \\
\end{pmatrix}.
\end{equation*}
The characteristic polynomial of $Q$ is given by
\begin{equation*}
f(\lambda)=\lambda^2 -2(n+l_1 -1)\lambda +{(n+l_1 -1)}^2 -4.
\end{equation*}
Let $\lambda_1$ and $\lambda_2$, $\lambda_1 \geq \lambda_2$, be the zeroes of $f (\lambda)$. It is easy to see that $\lambda_1 =n+l_1 +1$ and $\lambda_2 =n+l_1 -3$. Using Theorem \ref{T5}, we get $\partial^L_1 (G)\geq \lambda_1 =n+l_1 +1$. Also, since all $l_i$`s, $i=1,2,\dots,k$, are equal, that is, $l_1 =l_2 =\dots=l_k$ and  $n \equiv 0 \pmod k$, therefore, $l_1 =\frac{n}{k}  =\Big \lfloor \frac{n}{k}\Big \rfloor$. Thus,
\begin{equation*}
\partial^L_1 (G)\geq n+l_1 +1=n+\Big \lfloor \frac{n}{k}\Big \rfloor +1>n+\Big \lfloor \frac{n}{k}\Big \rfloor
\end{equation*}
Also, by Lemma \ref{L1}, $\partial^L_{n-1} (G)=n$. Hence,
\begin{equation*}
DLS(G)=\partial^L_1 (G)-\partial^L_{n-1} (G)>n+\Big \lfloor \frac{n}{k}\Big \rfloor-n=\Big \lfloor \frac{n}{k}\Big \rfloor.
\end{equation*}
This proves the result. \qed

\begin{corollary}
Let $G$ be a complete bipartite graph with $n$ vertices,  $n$ being an odd integer, then
\begin{equation*}
DLS(G)>\Big \lfloor \frac{n}{2}\Big \rfloor.
\end{equation*}
\end{corollary}
A graph $G$ is said to be $k-colorable$ if its vertex set $V$ can be colored with $k$ colors such that no two adjacent vertices get the same color. The \textit{chromatic number} of a graph $G$, denoted by $\chi$, is the minimum number of colors that can be used to color the vertices of $G$.

\begin{theorem}\label{T4}
\emph{ \cite{R7}} Let $G$ be a connected graph with $n\geq 4$ vertices and having chromatic number $\chi$ satisfying $\frac{n}{2}\leq \chi \leq n-1$. Then $\partial^L_1 (G)\geq \partial^L_1 (K_{\underbrace{{2,2,\dots,2}}_{n-\chi},\underbrace{{1,1,\dots,1}}_{2\chi -n}})$ with equality if and only if $G\cong K_{\underbrace{{2,2,\dots,2}}_{n-\chi},\underbrace{{1,1,\dots,1}}_{2\chi -n}}$.
\end{theorem}

Now, we show that $DLS(G)\geq 2$, for a connected graph $G$ with order $n\geq 4$ and given chromatic number.

\begin{theorem}\label{TA}
Let $G$ be a connected graph with $n\geq 4$ vertices having chromatic number $\chi$ satisfying $\frac{n}{2}\leq \chi \leq n-1$ and let complement of $G$ be disconnected. Then $DLS(G)\geq 2$ with equality if and only if $G\cong K_{\underbrace{{2,2,\dots,2}}_{n-\chi},\underbrace{{1,1,\dots,1}}_{2\chi -n}}$.
\end{theorem}
\noindent{\bf Proof.} Using Theorem \ref{T4} and noting that  $\partial^L_1 (K_{\underbrace{{2,2,\dots,2}}_{n-\chi},\underbrace{{1,1,\dots,1}}_{2\chi -n}})= n+2$, we observe that $\partial^L_1 (G)\geq n+2$ with equality if and only if $G\cong K_{\underbrace{{2,2,\dots,2}}_{n-\chi},\underbrace{{1,1,\dots,1}}_{2\chi -n}}$. As $\overline{G}$ is disconnected, using Lemma \ref{L1}, we get $\partial^L_{n-1}(G)=n$. Thus,
\begin{equation*}
DLS(G)=\partial^L_1 (G)-\partial^L_{n-1 }(G)\geq n+ 2-n=2
\end{equation*}
with equality  if and only if $G\cong K_{\underbrace{{2,2,\dots,2}}_{n-\chi},\underbrace{{1,1,\dots,1}}_{2\chi -n}}$. \qed

\noindent{\bf Acknowledgements.}  The research of S. Pirzada is supported by the SERB-DST research project number MTR/2017/000084.\\

\noindent{\bf Data availibility} Data sharing is not applicable to this article as no data sets were generated or analyzed
during the current study.\\	

\noindent{\bf Conflict of interest.} The authors declare that they have no conflict of interest.

\end{document}